\documentclass{amsart} %

\usepackage[mathscr]{eucal}
\usepackage{amsfonts,amssymb,amsmath,amsthm,graphicx}

\newtheorem{prop}{Proposition}
\newtheorem{cor}{Corollary}
\newtheorem{lem}{Lemma}
\newtheorem*{thm*}{Theorem}
\newtheorem*{cor*}{Corollary}

\theoremstyle{definition}

\newtheorem*{definition*}{Definition}

\newtheorem*{definitions*}{Definitions}

\newtheorem*{hist}{Historical note}
\newtheorem*{notation*}{Notation}
\newtheorem*{construction}{Construction} 
\newtheorem*{construction*}{Construction}
\newtheorem*{convention}{Convention}

\newcounter{xampl}
\newtheorem{example}[xampl]{Example}

     \def\a{\alpha}\def\b{\beta} \def\e{\varepsilon} \def\g{\gamma} 
       \def\phi{\varphi} 
     \def\s{\sigma} \def\t{\tau}
     \def\x{\xi}
     \def\D{\Delta} 
     
    \def\overstrike#1#2{{\setbox0\hbox{$#2$}\hbox to \wd0{\hss
                         $#1$\hss}\kern-\wd0\box0}}
\def\twobars#1#2#3#4#5#6{\vcenter{\hrule height.#1pt width#2pt
                               \vskip#3pt
                               \hrule height.#4pt width#5pt
                               \vskip#6pt}}
     \def\stroke#1#2#3{\vrule height#1pt width.#2pt depth#3pt}
     \def\connsum{\raise.25ex\hbox{\overstrike\parallel=}}
     \def\bdconnsum{\hskip2pt
                   \stroke83{-1.4}\twobars433431\stroke{5.1}31
                   \hskip2pt}
     
     \font\bxtwelve=cmbx12      
     \def\plumb#1{\thinspace{{\lower.75ex%
          \hbox{\text{\bxtwelve*}}}}_{#1}\thinspace} 
     \def\arg{\operatorname{arg}}
      
     \def\Bd{\partial} 
     \def\bword{\mathbf{b}} 
     \def\bydefl{\negthinspace:=\negthinspace}
     \def\bydefr{\negthinspace=:\negthinspace}
     \def\card#1{\operatorname{card}(#1)} 
     \def\C{\mathbb C} 
     \def\Cext{\P_1(\C)} 
     \def\crit{\operatorname{crit}} 
     \def\cut{\wr} 
     \def\emptyset{\varnothing}
     \def\from{\colon\thinspace} 
     \def\h#1#2{h^{\scriptscriptstyle(#1)}_{#2}} 
     
     \def\Im{\operatorname{Im}} 
     \def\Int{\operatorname{Int}} 
     \def\Nb#1#2{N_{#1}({#2})} 
     \def\MN{\mathscr M\mathscr N}
     \def\N{\mathbb N} 
     \def\P{\mathbb P} 
     \def\pr{\operatorname{pr}}

     \def\R{\mathbb R} 
     \def\rank{\operatorname{rank}}
     \def\Re{\operatorname{Re}} 
     \def\sub{\subset} 
     \def\Ts{\mathscr T}
	 \def\Z{\mathbb Z} 

\begin{document}

\thanks{Research partially supported by the Fonds National Suisse and by NSF grant DMS-9504832.}

\title[Murasugi sums of Morse maps to the circle]
{Murasugi sums of Morse maps to the circle, \\
Morse-Novikov numbers, and free genus of knots}
\author{Lee Rudolph}
\address{Department of Mathematics, Clark University, 
Worcester MA 01610 USA}  
\email{lrudolph@black.clarku.edu}

\begin{abstract} 
Murasugi sums can be defined as readily 
for Morse maps to $S^1$ of 
(arbitrary) link complements in $S^3$
as for fibrations over $S^1$
of (fibered) link complements in $S^3$.
As one application, I show that 
if a knot $K$ has free genus $m$, then 
there is a Morse map $S^3\setminus K\to S^1$ 
(representing the relative homology class 
of a Seifert surface for $K$)
with no more than $4m$ critical points.
\end{abstract}



\maketitle

\section{Introduction; statement of results}\label{introduction}

An oriented link $L\sub S^3$ determines a cohomology class 
$\xi_L\in H^1( S^3\setminus L;\Z)\cong %
\pi_0(\operatorname{Map}(S^3 \setminus L,K(\Z,1)))$.
The homotopy class of maps $S^3\setminus L\to S^1=K(\Z,1)$ 
corresponding to $\xi_L$ contains smooth maps which are \emph{Morse} 
(that is, all their critical points are non-degenerate),
and which restrict to a standard fibration in a 
neighborhood of $L$ (so they have only
finitely many critical points).  The \emph{Morse-Novikov number} 
$\MN(L)$ is the minimal number of critical points 
of such a map.

Tautologously, $\MN(L)=0$ 
iff $L$ is a fibered link.  It is natural to ask 
how $\MN(L)$ may be calculated---or estimated---for general $L$.
Moreover, for the class of fibered links, there exist both 
nice characterizations in other terms (e.g., 
a knot $K$ is fibered iff 
$\ker(\pi_1(S^3\setminus K)\to H_1(S^3\setminus K;\Z))$ 
is finitely generated)
and an array of interesting constructions (e.g.,
links of singularities,
Murasugi sums of fibered links,  
and closed homogeneous braids).  
Again, it is natural to ask what happens 
in general.

Some progress on these questions was made in \cite{P-R-W}. 
There, the Morse-Novikov theory 
of maps from manifolds to the circle (introduced by 
Novikov \cite{Novikov:multivalued}, and previously 
applied to knot complements in $S^3$ by Lazarev \cite{Lazarev}), 
and in particular the Novikov inequalities (analogues 
for Novikov homology of the classic Morse inequalities
for ordinary homology), were applied to give lower bounds for
$\MN(L)$: for example, it was shown that for any $n$, there 
exists a knot $K$ with $\MN(K)\ge n$; on the other hand,
it was also shown that there are knots for which the Novikov 
homology vanishes but the Morse-Novikov number does not.
An explicit construction established subadditivity of 
Morse-Novikov number over connected sum,
\begin{equation}\label{subadditivity over connected sum}
\MN(L_0\connsum L_1)\leq \MN(L_0)+\MN(L_1). \tag{$*$}
\end{equation}

Here, by constructing \emph{Murasugi sums of Morse maps}, 
I confirm the conjecture of \cite{P-R-W} that
a restatement of 
(\ref{subadditivity over connected sum})
in terms of Seifert surfaces 
extends to arbitrary Murasugi sums.
Applications include the following.  
(1)~For any knot $K$, $\MN(K)\le 4g_f(K)$, 
where $g_f(K)$ is the free genus of $K$;
for many knots $K$ (e.g., all but three
twist knots), $\MN(K)= 2g_f(K)$. 
(2)~The ``inhomogeneity'' of any braidword 
(or $\Ts$--bandword) with closed braid $L$
yields an upper bound on $\MN(L)$. 

Section~\ref{preliminaries} assembles preliminary material.
Section~\ref{Morse maps from Milnor maps} describes Milnor maps
and constructs a simple, but fundamental, example of a 
Morse map which is not a fibration.  
Section~\ref{Murasugi sums} constructs Murasugi sums 
of Morse maps and proves a generalization of 
(\ref{subadditivity over connected sum}).
Upper bounds on $\MN(L)$, 
and some exact calculations of $\MN(L)$ (in favorable cases),
are derived in sections~\ref{freeness} 
and \ref{inhomogeneity}.

Thanks to Andrei Pajitnov and Claude Weber
for helpful conversations.

\section{Preliminaries}\label{preliminaries}

The symbol $\square$ 
signals either the end or the omission of a proof.

Spaces, maps, etc., are smooth ($\mathscr C^\infty$) 
unless otherwise stated.
Manifolds may have boundary and are always oriented;
in particular, $\R$, $\C^n$, 
\[
D^{2n}\bydefl \{(z_1,\dots,z_n)\in\C^n:|z_1|^2+\dots+|z_n|^2\le 1\},
\]
and $S^{2n-1}\bydefl \Bd D^{2n}$ have standard orientations,
as does $S^2$ when it is identified with the Riemann sphere 
$\Cext\bydefl \C\cup\{\infty\}$. 
For suitable $Q\sub M$, let $\Nb MQ$ denote a closed 
regular neighborhood of $Q$ in $(M,\Bd M)$. 
If $Q$ is a codimension--$2$ submanifold of $M$ 
with trivial normal bundle, then a trivialization 
$\t\from Q\times D^2 \to \Nb{M}Q$ 
is \emph{adapted to} a map $f\from M\setminus Q\to S^1$ if
$\t(\x,0)=\x$ and $f(\t(\x,z)) = z/|z|$ 
for $\x\in Q$, $z \in D^2 \setminus \{0\}$.

A \emph{surface} is a compact 
$2$--manifold no component of which has empty boundary.  
A \emph{Seifert surface} is a surface $F\sub S^3$.  
A \emph{link} $L$ is the boundary of a Seifert 
surface $F$, and any $F$ with $\Bd F=L$ is a
Seifert surface \emph{for} $L$.  
A \emph{knot} is a connected link.

Let $L\subset S^3$ be a link.
The image in $H_2(S^3,L;\Z)$ 
of the fundamental class $[F] \in H_2(F,L;\Z)$ 
is independent of the choice of a Seifert surface $F$ for $L$;
denote the corresponding orientation class for $L$
in $H^1(S^3\setminus L;\Z)\cong H_2(S^3,L;\Z)$
by $\xi_L$.
Call $f\from S^3 \setminus L \to S^1$ \emph{simple}
if $[f] \cong \xi_L$ in 
$\pi_0(\operatorname{Map}(S^3 \setminus L,S^1))\cong %
H^1(S^3 \setminus L; \Z)$.
The \emph{Morse-Novikov number of $L$} is 
$\MN(L)\bydefl \min\{n: \textrm{there exists a simple Morse}$ %
$\textrm{map } f\from S^3 \setminus L \to S^1 %
\textrm{ with exactly $n$ critical points}\}$.  On general 
principles, $\MN(L)$ $<\infty$.

\begin{definitions*}\label{properties of Morse maps}
A simple Morse map $f\from S^3 \setminus L \to S^1$ is 
\begin{enumerate}

\item \label{boundary-regular}
\emph{boundary-regular} if there is a trivialization 
$\t\from L\times D^2 \to \Nb{S^3}L$ 
which is adapted to $f$; 

\item \emph{moderate} if no critical point of $f$ is 
a local extremum;

\item \emph{self-indexed} if the value of $f$ at a critical point 
is an injective function of the index of the critical point; 

\item 
\emph{minimal} if no Morse map in $\xi_L$
has strictly fewer critical points than $f$.

\end{enumerate}

\end{definitions*}

\begin{prop}\label{implications among properties of Morse maps}
Let $f\from S^3 \setminus L \to S^1$ be a simple 
Morse map.  
\begin{enumerate}
\item\label{finitely many cps is boundary-regular} 
If $f$ has only finitely many critical points, then 
$f$ is properly isotopic to a boundary-regular simple Morse map. 
\item\label{boundary-regular gives Seifert surfaces}  
If $f$ is boundary-regular then 
$f$ has only finitely many critical points,
and for every regular value $\exp(i\theta)\in S^1$, 
$S(f,\theta)\bydefl L\cup f^{-1}(\exp(i\theta))$ 
is a Seifert surface for $L$.
\item\label{minimal is moderate}
If $f$ is minimal, then $f$ is moderate
and has precisely $\MN(L)$ critical points.  
\item\label{moderate is self-indexed} 
If $f$ \emph{(}has only finitely many critical points
and\emph{)} is moderate, then up to \emph{(}proper\emph{)} isotopy 
$f$ is moderate and self-indexed.
\item\label{self-indexed} If $f$ is boundary-regular,
moderate, 
and self-indexed,  
then either
\begin{enumerate}
\item 
\label{fibration}
$f$ has no critical points
\emph{(\emph{so $\MN(L)=0$})} and is in fact a fibration, 
and the Seifert surfaces $S(f,\theta)$ 
are mutually isotopic in $S^3$, or
\item 
\label{non-fibration} $f$ has $2m>0$ critical points,
$m$ each of index $1$ and index $2$,
and the Seifert surfaces $S(f,\theta)$ 
fall into exactly two distinct isotopy classes in $S^3$.
\end{enumerate}
\end{enumerate}
\end{prop}

\begin{proof} 
Mostly straightforward (for (\ref{minimal is moderate}) 
and (\ref{non-fibration}),
see \cite{P-R-W}).
\end{proof}

In case (\ref{fibration}), 
$L$ is (as usual) called 
a \emph{fibered link} and $S(f,\theta)$ is called 
a \emph{fiber surface of $f$}.  (Also as usual, 
any Seifert surface for $L$ isotopic to an 
$S(f,\theta)$ is called a \emph{fiber surface for $L$}.)
In case (\ref{non-fibration}), the two isotopy classes of 
Seifert surfaces $S(f,\theta)$ for $L$ (which could well be 
a fibered link even though $f$ is not a fibration) are 
distinguished by their Euler characteristics,
which differ by $2m$.
Any of the Seifert surfaces $S(f,\theta)$ 
with first homology group of larger (resp., smaller) rank
will be called a \emph{large} 
(resp., a \emph{small}) \emph{Seifert surface of} 
the boundary-regular moderate self-indexed simple Morse map $f$ 
(note: not ``of $L$'').
Any fiber surface of a fibration $f$ may be called
either large or small.

\begin{convention}
Henceforth, all Morse maps are boundary-regular and simple.
\end{convention}

\section{Morse maps from Milnor maps}\label{Morse maps from Milnor maps}

The first explicit Morse maps (in fact, fibrations)
for an infinite class of links were given by 
Milnor's celebrated Fibration Theorem \cite{Milnor:singular-points}, 
where they appear as (the instances for functions $\C^2\to \C$ of) what 
are now called the ``Milnor maps'' associated
to singular points of complex analytic functions $\C^n\to\C$.  
For present and future purposes, it is useful to somewhat
extend the framework in which Milnor studied these maps.
Given a non-constant meromorphic function $F\from M\dashrightarrow \Cext$
on a complex manifold $M$, let $D(F)$ be the (possibly singular)
complex hypersurface which is the closure in $M$ of 
$F^{-1}(0)\cup F^{-1}(\infty)$. 
The \emph{argument} of $F$ is 
$\arg(F)\bydefl F/|F|\from M\setminus D(F) \to S^1$.
For $M=\C^n$, the restriction of $\arg(F)$ to $rS^{2n-1}\setminus D(F)$
is the \emph{Milnor map of $F$ at radius $r$}.
Call the Milnor map of $F$ at radius $1$ simply 
\emph{the Milnor map of $F$}.

Milnor's proof of  \cite[Lemma~4.1]{Milnor:singular-points},
stated for analytic $F$, applies equally well to meromorphic
$F$. 

\begin{lem}\label{Milnor's lemma}
If $F\from \C^n\dashrightarrow\Cext$ is meromorphic,
then $(z_1,\dots,z_n)\in rS^{2n-1}\setminus D(F)$
is a critical point of the Milnor map of $F$ at radius $r$
iff the complex vectors 
\[
(\overline{z}_1,\dots,\overline{z}_n), \medspace
\frac{1}{i\thinspace F(z_1,\dots,z_n)}%
\bigl(\frac{\partial F}{\partial z_1}(z_1,\dots,z_n),\dots,
\frac{\partial F}{\partial z_1}(z_1,\dots,z_n)\bigr)
\]
are linearly dependent over $\R$.
\qed
\end{lem}

The next lemma rephrases several more results from 
\cite{Milnor:singular-points}.

\begin{lem}\label{Milnor maps are nice}
Let $F\from \C^2\to\C$ be analytic
\emph{(\emph{so $D(F)=F^{-1}(0)$})}, with $F(0,0)=0$,
and suppose $F$ has no repeated factors in ${\mathscr O}_{\C^2}$
\emph{(\emph{so the multiplicity of each irreducible
component of $D(F)$ is $1$})}.
\begin{enumerate}
\item\label{most radii give a link} 
There is a set $X(F)$ of radii $r\in{]}0,\infty{[}$ 
such that $X(F)$ is discrete in $[0,\infty{[}$ and,
if $r\notin X(F)$, then $D(F)$ intersects $rS^3$ 
transversally, so that $L(F,r)\bydefl (1/r)(D(F)\cap rS^3)$ 
is a link in $S^3$; $L(F,r)$ and $L(F,r')$ are isotopic
for $r, r'$ in the same component of ${]}0,\infty{[}\setminus X(F)$.
\item\label{Milnor maps are simple and adaptable} 
If $1\notin X(F)$, then the Milnor map of $F$ is 
simple, and $\Nb{S^3}{L(F,1)}$ has a trivialization 
which is adapted to the Milnor map of $F$.
\item\label{Milnor's fibration theorem for curves}
If $1<\inf X(F)$, then the Milnor map of $F$ is
a fibration \emph{(\emph{the} Milnor fibration of $F$)},
and in particular $L(F,1)$ \emph{(\emph{the} link of the 
singularity of $F^{-1}(0)$ at $(0,0)$)} is a fibered link.
\qed
\end{enumerate}
\end{lem}

Note that, in case (\ref{Milnor maps are simple and adaptable}),
the Milnor map of $F$ may have degenerate critical
points, but if it does not, then it is Morse (in the sense of
the convention in section~\ref{preliminaries}).  
In practice, Lemma~\ref{Milnor's lemma}
makes it easy to locate the critical points of a Milnor map
and check them for non-degeneracy.

\begin{example}\label{o{m,n} as a Milnor fibration}
Given $m,n\ge 0$ with $m+n=1$ in case $mn=0$,
let $F_{m,n}\from \C^2\to\C : (z,w)\mapsto mz^m+nw^n$; easily,
$X(F_{m,n})=\emptyset$. 
The link of the singularity of $F_{0,1}^{-1}(0)$ at $(0,0)$
is an \emph{unknot} $O$ (of course 
$F_{0,1}^{-1}(0)$ is not in fact singular at $(0,0)$); 
the Milnor fibration of $F_{0,1}=\pr_2$ will be denoted by $o$.
More generally, the link of the 
singularity of $F_{m,n}^{-1}(0)$ at $(0,0)$ is a \emph{torus link 
$O\{m,n\}$ of type $(m,n)$}, and the Milnor fibration of $F_{m,n}$ 
will be denoted by $o\{m,n\}$.  
\end{example}

\begin{figure}
\centering
\includegraphics[width=3.5in]{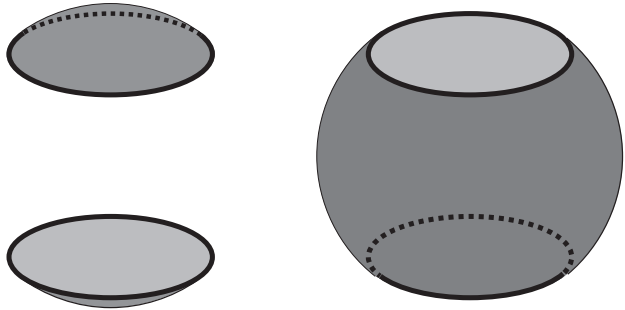}
\caption{A small and a large Seifert surface of $u$.
\label{uround}}
\end{figure}

\begin{example}\label{u as a Milnor Morse map}
If $G_\e\from \C^2\to\C : (z,w)\mapsto 4w^2-8\e w-1$, then 
for sufficiently small $\e\ge 0$, $U_\e\bydefl L(G_\e,1)$ is an unlink
of two unknotted components, and if $\e> 0$, then 
(by direct calculation using Lemma~\ref{Milnor's lemma}) 
the Milnor map of $G_\e$ is Morse, 
with two critical points (one each of index $1$ and index $2$).
Let $U\bydefl U_\e$ for a sufficiently small $\e>0$, and denote
by $u\from S^3\setminus U\to S^1$ this Morse Milnor map.  (The situation
is pictured in Fig.~\ref{uround}.)  
\end{example}

\begin{prop}\label{MN(U)=2}
$\MN(U)=2$ and $u$ is a minimal Morse map.
\end{prop}

\begin{proof}
This is immediate from the properties of $u$ just asserted,
given that $U$ is not fibered (which follows, for instance,
from the fact that $\pi_2(S^3\setminus U)\ne \{0\}$).
\end{proof}

The non-fibration $u$, trivial though it be, is an ingredient
of fundamental importance in the constructions of section~\ref{freeness}.

\section{Murasugi sums of Morse maps}\label{Murasugi sums}

For $n\in\N$, write $G_{n}\bydefl \{z\in\C:z^{n}=1\}$, and 
let $p_n\from \Cext\setminus G_{2n}\to S^1$
be the argument of the rational function
$\Cext\to\Cext : z\mapsto(1+z^n)/(1-z^n)$, viz., 
\[
p_n(z)=\frac{(1+z^n)/|1+z^n|}{(1-z^n)/|1-z^n|} 
\textrm{ for $z\notin G_{2n}\cup\{\infty\}$, }
p_n(\infty)=-1.
\]
Define $P_n\from \Cext\setminus G_{2n}\times[0,\pi]\to S^1$ 
by $P_n(z,\theta)=\exp(i\theta)p_n(z)$.  Note that
$p_n(\zeta z)=p_n(z)$ and $P_n(z,\theta)=P_n(\zeta z,\theta)$
for any $\zeta\in G_n$.

\begin{lem}\label{description of pn and Pn}
\emph{(1)}~$p_1$ is a fibration with fiber 
${]}-1,1{[}=p_1^{-1}(1)$.
For $n>1$, $p_n$ has exactly two critical points, 
$0\in p_n^{-1}(1)$ and $\infty\in p_n^{-1}(-1)$,
at each of which the germ of $p_n$ is smoothly 
conjugate to the germ of $z\mapsto\Im(z^n)$ at $0$.
\emph{(2)}~$P_1$ is a fibration with fiber 
${]}-1,1{[}\times[0,\pi] =P_1^{-1}(1)$.
For all $n$, $P_n$ has no critical points,
and there is a trivialization 
$\t_{P_n}\from (G_{2n}\times[0,\pi])\times D^2%
\to\Nb{\Cext\times[0,\pi]}{G_{2n}\times[0,\pi]}$ 
which is adapted to $P_n$.
\qed
\end{lem}

\begin{figure}
\centering
\includegraphics[width=3.5in]{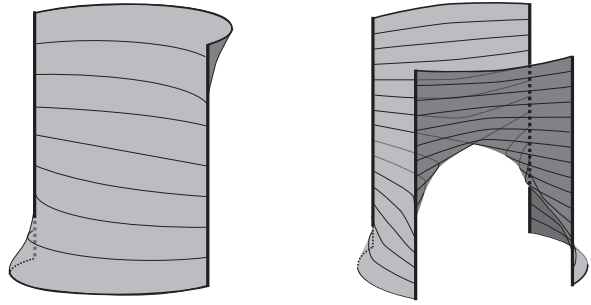}
\caption{$Q_1(\pi/2)$ and some 
level sets of $\pr_2|Q_1(\pi/2)$; 
$Q_2(\pi/2)$ and some level sets of $\pr_2|Q_2(\pi/2)$.
\label{q1+q2}}
\end{figure}

Let $Q_n(\theta)$ denote the closure of $P_n^{-1}(\exp(i\theta))$ 
in $\Cext\times[0,\pi]$.  
By inspection (and Lemma~\ref{description of pn and Pn}), 
for all $\exp(i\theta)\in S^1$, $Q_1(\theta)$ is a $4$--gonal $2$--disk
with smooth interior bounded by the union of 
$G_{2}\times[0,\pi]$ and $2$ semicircles in $\Cext\times\{0,\pi\}$;
for $n>1$ and all $\exp(i\theta)\in S^1\setminus\{1,-1\}$, 
$Q_n(\theta)$ is a $4n$--gonal $2$--disk
with smooth interior bounded by the union of
$G_{2n}\times[0,\pi]$ and $2n$ circular arcs 
in $\Cext\times\{0,\pi\}$.  
(The cases $n=1, 2$, $\theta=\pi/2$, are pictured in 
Fig.~\ref{q1+q2}.)

\begin{lem}
\label{description of Qn}
For all $\exp(i\theta)\in S^1$,
the restriction $\pr_2|Q_1(\theta)\from Q_1(\theta)\to[0,\pi]$ has 
no critical points.  
For $n>1$ and $\exp(i\theta)\in S^1\setminus\{1,-1\}$, 
there is exactly one critical point 
of $\pr_2|Q_n(\theta)\from Q_n(\theta)\to[0,\pi]$ 
\emph{(}to wit,
$(0,\theta)$ for $0<\theta<\pi$
and $(\infty,\theta-\pi)$ for $\pi<\theta<2\pi$\emph{)},
at which the germ 
of $\pr_2|Q_n(\theta)$ is smoothly conjugate to the germ of 
$z\mapsto\Im(z^n)$ at $0$. \qed
\end{lem}

\begin{figure}
\centering
\includegraphics[width=3.5in]{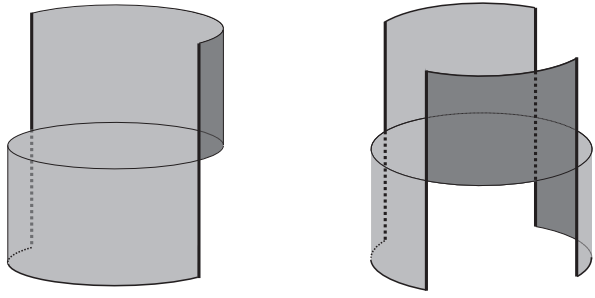}
\caption{$Q'_1(\pi/2)$ and $Q'_2(\pi/2)$.\label{q1'&q2'}}
\end{figure}

\begin{figure}
\centering
\includegraphics[width=3.5in]{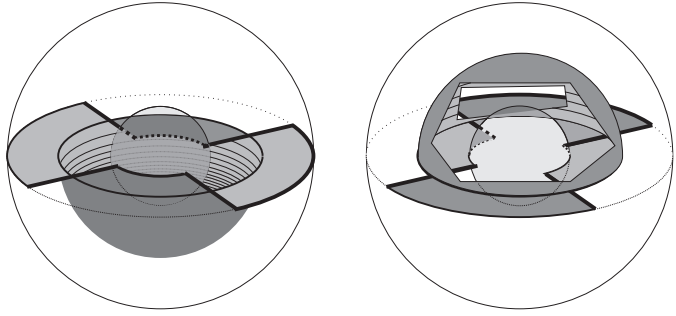}
\caption{$Q'_2(\pi/2)$ and $Q'_2(3\pi/2)$ 
(the latter with viewports),
as they appear in in 
$\Nb{{S^3}}{{S^2}}\protect\cong \widehat\C\times[0,\pi]$.
\label{spherical q2'}}
\end{figure}

By further inspection, for
$0<\theta<2\pi$, $\theta\ne\pi$, 
$Q_n(\theta)$ is piecewise-smoothly isotopic
(by an isotopy fixing 
$\Bd Q_n(\theta)\cup (Q_n\cap %
\Cext\times\{\theta-\lfloor \theta/\pi \rfloor\pi\})$
pointwise) to a piecewise-smooth \hbox{$4n$--gonal}
\hbox{$2$--disk} $Q'_n(\theta)$ 
which is so situated 
that both
$Q'_n(\theta)\cap 
(\Cext\times %
[0,\theta-\lfloor \theta/\pi \rfloor\pi])$ 
and 
$Q'_n(\theta)\cap %
(\Cext\times %
[\theta-\lfloor \theta/\pi \rfloor\pi,\pi])$ 
are piecewise-smooth {$4n$--gonal} $2$--disks,
while $Q'_n(\theta)\cap %
(\Cext\times \{\theta-\lfloor \theta/\pi \rfloor\pi\})$ 
is a (smooth) $2$--disk in 
$\Cext\times\{\theta-\lfloor \theta/\pi \rfloor\pi\}$ 
naturally endowed with the structure of a \hbox{$2n$--gon}.
(The cases $n=1, 2$, $\theta=\pi/2$, are pictured in 
Fig.~\ref{q1'&q2'}.
Versions of $Q'_2(\pi/2)$ and $Q'_2(3\pi/2)$
in $\Nb{{S^3}}{{S^2}}\cong \Cext\times[0,\pi]$ 
are pictured in Fig.~\ref{spherical q2'}.)

\begin{definitions*}\label{f-good neighborhoods}
Let $L \sub S^3$ be a link, 
$f\from S^3\setminus L\to S^1$ a 
Morse map, 
$\exp(i\theta)\in S^1$ a regular value of $f$,
and $\psi\sub S(f,\theta)$ an \emph{$n$--star} on 
$S(f,\theta)$ (in the sense of \cite{constqp5}:
that is, $\psi$ is the union of $n$ arcs $\a_s$,
pairwise disjoint except for a common endpoint
$*_\psi\in\Int S(f,\theta)$, such that 
$\Bd S(f,\theta)\cap\a_s=\Bd\a_s\setminus\{*_\psi\}$
for each $s$).
A regular neighborhood $\Nb{S^3}\psi$ is \emph{$f$--good}
provided that $\Nb{S^3}\psi \cap S(f,\theta)$ is a regular neighborhood
$\Nb{S(f,\theta)}\psi$ (and thus an \emph{$n$--patch} on 
$S(f,\theta)$ in the sense of \cite{constqp5}:
that is, a $2$--disk naturally endowed with 
the structure of a \hbox{$2n$--gon} whose edges are 
alternately boundary arcs and proper arcs in $S(f,\theta)$)
and there is a diffeomorphism 
$h\from (\Cext,G_{2n}) \to (\Bd\Nb{S^3}\psi, L \cap \Bd\Nb{S^3}\psi)$
such that $(f\circ h)|(\Cext\setminus G_{2n})=p_n$.
\end{definitions*}

\begin{lem}\label{good patches exist}
Every neighborhood of $\psi$ in $S^3$ contains an 
$f$--good regular neighborhood.
\end{lem}
\begin{proof} 
First suppose that $L=O$ and $f=o$ (as in 
Example~\ref{o{m,n} as a Milnor fibration}).
Stereographic projection 
$\s\from S^3\setminus\{(0,-i)\} \to\C\times\R 
        : (z,w) \mapsto(z,\Re(w))/(1+\Im(w))$ 
maps $O$ to $S^1\times\{0\}$ and $S(o,\pi/2)$ to
$D^2\times\{0\}$. The radius 
\hbox{$\psi_1\bydefl \s^{-1}([0,1]\times\{0\})$}
is a \hbox{$1$--star} in $S(o,\pi/2)$, and 
the preimage $\s^{-1}(B)$ of an appropriate ellipsoidal 
\hbox{$3$--disk} $B\sub\C\times\R$ (say, with one focus at $(0,0)$,
center at $(1-\e,0)$ for sufficiently small $\e>0$,
and minor axes much shorter than the major axis),
is an $o$--good regular neighborhood $\Nb{S^3}{\psi_1}$.
The $n$--sheeted cyclic branched covering 
\[
(\C\times\R)\cup\{\infty\}\to(\C\times\R)\cup\{\infty\} %
: (\zeta,t)\mapsto(\zeta^n,t),
\infty\mapsto\infty,
\]
branched along $(\{0\}\times\R)\cup\{\infty\}$, 
can be modified in a neighborhood of $\infty$ 
so as to induce via $\s^{-1}$ a smooth 
$n$--sheeted cyclic branched covering $c_n\from S^3\to S^3$,
branched along $A\bydefl \{0\}\times S^1\sub S^3$, with 
$c_n^{-1}(S(o,\pi/2))=S(o,\pi/2)$; then
\hbox{$\psi_n\bydefl c_n^{-1}(\psi_1)$} is an $n$--star
in $S(o,\pi/2)$, and $c_n^{-1}(\Nb{S^3}{\psi_1})$
is (if minimal care has been taken) an 
$o$--good regular neighborhood $\Nb{S^3}{\psi_n}$.
(The cases $n=1,2,3$ are pictured in 
Fig.~\ref{o-good neighborhoods}.) 
Clearly any neighborhood of $\psi_n$ in $S^3$
contains an $o$--good regular neighborhood of the type
just constructed.

\begin{figure}
\centering
\includegraphics[width=4.25in]{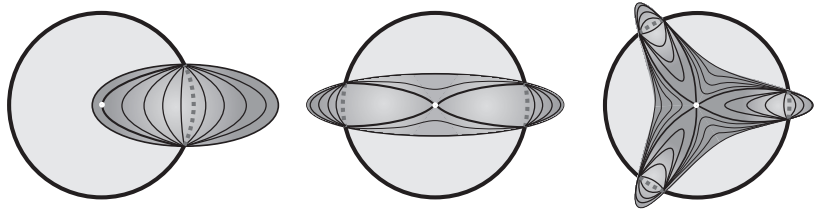}
\caption{Some level sets of $o|(\Bd\Nb{S^3}{\psi_n}\setminus O)$
for $n=1,2,3$.
\label{o-good neighborhoods}}
\end{figure}

The general case follows immediately, upon observing that,
for any link $L$, 
Morse map $f\from S^3\setminus L\to S^1$, 
and regular value $\exp(i\theta)$ of $f$, 
if $\psi\sub S(f,\theta)$ is an \hbox{$n$--star} on $S(f,\theta)$, 
then there is a diffeomorphism $h\from M\to h(M)$ 
from a (non--$f$--good) neighborhood $M$ of $\psi$ in $S^3$ 
to a neighborhood $h(M)$ of $\psi_n$ in $S^3$,
and a diffeomorphism $k\from (S^1,\exp(i\theta))\to (S^1,i)$,
such that $h(M\cap S(f,\theta))=h(M)\cap S(o,\pi/2)$
and $k\circ f|(M\setminus L)=o\circ h|(M\setminus L)$.
(A ``side view'' of part of one such $M$ is pictured in 
Fig.~\ref{sideview}.)  
\begin{figure}
\centering
\includegraphics[width=2in]{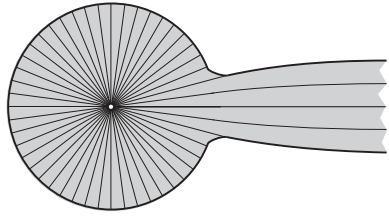}
\caption{A ``side view'' (indicating some level sets
of $f$) of the part of the non--$f$--good neighborhood 
$M$ near an endpoint of $\psi$.\label{sideview}}
\end{figure}
\end{proof}

\begin{construction}\label{M. sum of Morse maps}
For $s=0,1$,
let $L_s \sub S^3_s$ be a link and 
\hbox{$f_s\from S^3_s\setminus L_s\to S^1$} a Morse map
with critical points $\crit(f_s)\sub S^3_s\setminus L_s$.
Fix $n>0$ and
$\exp(i\Theta)\in %
S^1\setminus(f_0(\crit(f_0))\cup f_1(\crit(f_1))\cup\{0,\pi\})$.
Let $\psi_{n,s}\sub S(f_s,\Theta)$ be an $n$--star
with $f_s$--good regular neighborhood $N_s\bydefl \Nb{S^3_s}{\psi_{n,s}}$.
Let $E_s\bydefl S^3_s\setminus\Int N_s$; let
$d_s\from E_s\to D^3$ be a diffeomorphism.
Let $h_s\from (\Cext,G_{2n}) \to (\Bd N_s, L \cap \Bd N_s)$
be a diffeomorphism with 
$(f_s\circ h_s)|(\Cext\setminus G_{2n})=p_n$.
For a fixed $\zeta\in G_n$, let 
$(z,s)\equiv_\zeta h_s(\zeta^s z,s)$ 
($z\in \Cext, s=0,1$).  

The identification space $\Sigma(N_0,N_1,\zeta)\bydefl  %
(E_0 \sqcup %
\thinspace \Cext\times[0,\pi] \thinspace %
\sqcup E_1 ){/}{\equiv}_\zeta$
has a natural piecewise-smooth structure,
imposed on it by the identification map $\Pi$,
with respect to which the $1$--submanifold 
\[
L(N_0,N_1,\zeta)\bydefl (L_0\cap E_0 \sqcup %
G_{2n}\times[0,\pi]\sqcup L_1\cap E_1){/}{\equiv}_\zeta
\]
and the map 
\[
f(N_0,N_1,\zeta)\bydefl %
((f_0 \sqcup P_n \sqcup f_1){/}{\equiv}_\zeta)|
(\Sigma(N_0,N_1,\zeta)\setminus L(N_0,N_1,\zeta)) 
\]
are smooth where this is meaningful (i.e., in the complement
of the identification locus $\Pi(\Cext\times\{0,\pi\})$,
along which $\Sigma(N_0,N_1,\zeta)$ is itself 
\emph{a priori} only piecewise-smooth).
It is not difficult to give $\Sigma(N_0,N_1,\zeta)$ 
a smooth structure everywhere, in which 
$L(N_0,N_1,\zeta)$ and $f(N_0,N_1,\zeta)$ 
are everywhere smooth.
The smoothing can be done very naturally off
$\Pi((G_{2n}\cup\{0,\infty\})\times\{0,\pi\})$, 
using the fact that by construction 
(see Lemmas~\ref{description of pn and Pn} and
\ref{description of Qn}) 
the map $f(N_0,N_1,\zeta)$ 
and the function 
$(\|d_0\|^2-1) \sqcup \pr_2 \sqcup (\pi+1-\|d_1\|^2){/}{\equiv}_\zeta$
are piecewise-smoothly transverse (in an evident sense)
on $\Pi((\Cext\setminus(G_{2n}\cup\{0,\infty\}))\times\{0,\pi\})$.
A somewhat less natural, but not difficult,
construction smooths $\Sigma(N_0,N_1,\zeta)$ on 
$\Pi(\{0,\infty\}\times\{0,\pi\})$ in such a way 
as to make $f(N_0,N_1,\zeta)$ smooth there also
(in the process, possibly forcing
$(\|d_0\|^2-1) \sqcup \pr_2 \sqcup (\pi+1-\|d_1\|^2){/}{\equiv}_\zeta$
to be not smooth at those points).
Nor is there is any difficulty in 
smoothing $\Sigma(N_0,N_1,\zeta)$ 
on $\Pi(G_{2n}\times\{0,\pi\})$.  
Further details will be suppressed.

When $\Sigma(N_0,N_1,\zeta)$, with the smooth structure just
constructed, is identified with $S^3$, 
$L(N_0,N_1,\zeta)$ (resp., $f(N_0,N_1,\zeta)$) 
will be called a \emph{$2n$--gonal Murasugi sum} 
of $L_0$ and $L_1$ (resp., of $f_0$ and $f_1$) 
and denoted by $L_0\plumb{(N_0,N_1,\zeta)}L_1$ 
(resp., $f_0\plumb{(N_0,N_1,\zeta)}f_1$) or simply 
by $L_0\plumb{}L_1$ (resp., $f_0\plumb{}f_1$).
(A $2$--gonal Murasugi sum of links is simply a 
connected sum, and a $2$--gonal Murasugi sum of Morse
maps $f_0$ and $f_1$ may be denoted $f_0\connsum f_1$.
A more detailed description of the construction
in the case of connected sums, which goes more 
smoothly than the general case, is given in \cite{P-R-W}.)
\end{construction}

\begin{prop}\label{properties of Murasugi sum of maps}
If $L_s \sub S^3_s$ is a link and 
$f_s\from S^3_s\setminus L_s\to S^1$ is a Morse map
$(s=0,1)$, then 
\begin{enumerate}
\item\label{sum of links is a link}
any Murasugi sum $L_0\plumb{(N_0,N_1,\zeta)}L_1$ is a link, and 
\item\label{sum of maps is a map}
any Murasugi sum $f_0\plumb{(N_0,N_1,\zeta)}f_1$ is a Morse map
and its critical points 
$\crit(f_0\plumb{(N_0,N_1,\zeta)}f_1)=\crit(f_0)\cup\crit(f_1)$
have indices and critical values that are inherited unchanged
from those of $f_1$ and $f_2$;
in particular, 
\begin{enumerate}
\item
\label{sum preserves moderation}
if $f_0$ and $f_1$ are moderate, 
then $f_0\plumb{(N_0,N_1,\zeta)}f_1$ is moderate, and
\item 
\label{sum preserves self-indexing}
if $f_0$ and $f_1$ are self-indexed
and have the same critical value for each index,
then $f_0\plumb{(N_0,N_1,\zeta)}f_1$ is self-indexed.
\qed
\end{enumerate}
\end{enumerate}
\end{prop}

With $\Theta$ chosen as in the construction, the Seifert surface 
$S(f_0\plumb{(N_0,N_1,\zeta)}f_1,\Theta)$ 
is isotopic (up to smoothing) to 
$\Pi((S(f_0,\Theta)\cap E_0)\sqcup Q_n(\Theta) %
 \sqcup (S(f_1,\Theta)\cap E_1))$, and thus
piecewise-smoothly isotopic to
\begin{multline*}\label{}
\Pi(S(f_0,\Theta)\cap E_0 %
   \sqcup Q'_n(\Theta) %
   \sqcup S(f_1,\Theta)\cap E_1) = \\
                  \Pi(S(f_0,\Theta)\cap E_0\sqcup %
                  Q'_n(\Theta)\cap %
                  \Cext\times [0,\Theta-\lfloor \Theta/\pi %
                  \rfloor\pi]) \medspace\cup\qquad\quad\\
\Pi(Q'_n(\Theta)\cap %
                  \Cext\times [\Theta-\lfloor \Theta/\pi %
                  \rfloor\pi,\pi]
                  \sqcup S(f_1,\Theta)\cap E_1)\bydefr S'_0\cup S'_1
\end{multline*}
where $S'_s$ is piecewise-smoothly isotopic to $S(f_s,\Theta)$
by an isotopy carrying the $2$--disk $S'_0\cap S'_1$ (with
its naturally structure of $2n$--gon, noted after 
Lemma~\ref{description of Qn}) to the $n$--patch
$N_s\cap S(f_s,\Theta)$.
Thus, for a suitable 
diffeomorphism
$H\from N_0\cap S(f_0,\Theta)\to N_1\cap S(f_1,\Theta)$
(determined by $\zeta$, up to isotopy), 
$S(f_0\plumb{(N_0,N_1,\zeta)}f_1,\Theta)$ 
is a $2n$--gonal Murasugi sum 
$S(f_0,\Theta)\plumb{H}S(f_1,\Theta)$
as described in \cite{constqp5}
(see also the primary sources 
\cite{Murasugi:plumbing,Stallings,Gabai:Murasugi1}).

On the level of Seifert surfaces, 
a $2$--gonal Murasugi sum $S_0\plumb{H}S_1$
is the same as a boundary-connected sum 
$S_0\bdconnsum S_1$.  After boundary-connected sum, 
the most commonly encountered case of Murasugi sum is
$4$--gonal, and the most familiar and probably most
useful $4$--gonal Murasugi sums are \emph{annulus plumbings},
as described in \cite{espaliers} (see also the primary
sources \cite{Siebenmann:rational,Gabai:arborescent},
as well as \cite{Hayashi-Wada:plumbing} and
\cite{Sakuma:specialarborescent}).

Specifically, for any knot $K$, let $A(K,n)$ denote
any Seifert surface diffeomorphic to an annulus such that 
$K\subset\Bd A(K,n)$ and the Seifert matrix of $A(K,n)$
is $[n]$.  (For example, 
$A(O,-1)$ is a \emph{positive Hopf annulus}, isotopic
to a fiber surface of $o\{2,2\}$, as in 
Example~\ref{o{m,n} as a Milnor fibration}; 
the mirror image $A(O,1)$ of $A(O,-1)$
is a \emph{negative Hopf annulus}; and
$A(O,0)$ is isotopic to a large Seifert surface of $u$
as in Example~\ref{u as a Milnor Morse map}.)
Let $\g(K)\sub A(K,n)$ be an arc 
from $K$ to $\Bd A(K,n)\setminus K$,
so that $\Nb{A(K,n)}{\g(K)}\bydefr C_1$ is a $2$-patch (though
$\g(K)$, as given, is not a $2$-star).
Let $S$ be a Seifert surface $S$, $\a\sub S$ a proper arc,
$C_0\sub S$ a $2$--patch 
with $\a\sub\Bd C_0$ (respecting orientation).
Let $H\from (C_0,\a)\to (C_1,C_1\cap K)$ be a diffeomorphism.
Each of $\g(K)$, $C_0$, and $H$ is unique up to isotopy,
so the $4$--gonal Murasugi sum $S\plumb{H} A(K,n)$ 
depends (up to isotopy) only on $\a$, 
and there is no abuse in denoting it by 
$S\plumb{\a} A(K,n)$.  When $S=A(K',n')$ is also an 
annulus, it is slightly abusive---but handy---to write 
simply $A(K',n')\plumb{}A(K,n)$, with the understanding 
that $\a=\g(K')$.

\begin{example}\label{plumbing A(O,0)}
Given any proper arc $\a$ on a surface $S$,
there is a (highly non-unique!) $(0,1)$--handle decomposition 
\begin{equation}\label{handle decomposition}
S=\bigcup_{x\in X}\h0x\cup \bigcup_{z\in Z}\h1z\tag{$\dagger$}
\end{equation}
for which $\a$ is an attaching arc of some $1$--handle $\h1s$.
The subsurface 
\[
S\cut\a \bydefl \bigcup_{s\ne x\in X}\h0x\cup \bigcup_{z\in Z}\h1z
\]
is independent, up to isotopy on $S$, of the handle
decomposition (\ref{handle decomposition}), and 
can be thought of as ``$S$ cut open along $\a$''.
If $S$ is a Seifert surface, then for every $n\in\Z$ 
there is a reimbedding $T^n_\a\from S\to S^3$ (unique up
to isotopy) which is the identity on $(S\cut\a)\cup \kappa$,
where $\kappa\sub\h1s$ is a core arc, and which 
puts $n$ counterclockwise full twists in $\h1s$.
For any $n$, the annulus-plumbed Seifert surfaces
$S\plumb{\a} A(O,0)$ and $T^n_\a(S)\plumb{\a} A(O,0)$
are isotopic, and the boundary of each is isotopic to
$\Bd(S\cut\a$).
(The situation is pictured in 
Fig.~\ref{twist across alpha}, for $n=1$.)
In particular, for any knot $K$ and integer $n$,
$A(K,n)=T^n_{\g(K)}(A(K,0))$, 
and $\Bd (A(K,n)\plumb{}A(O,0))$ is an unknot.
\begin{figure}
\vspace{.5in}
\setlength{\unitlength}{1in}
\begin{picture}(4.5,2.25)(-.25,-.05)
\includegraphics[width=4.5in]{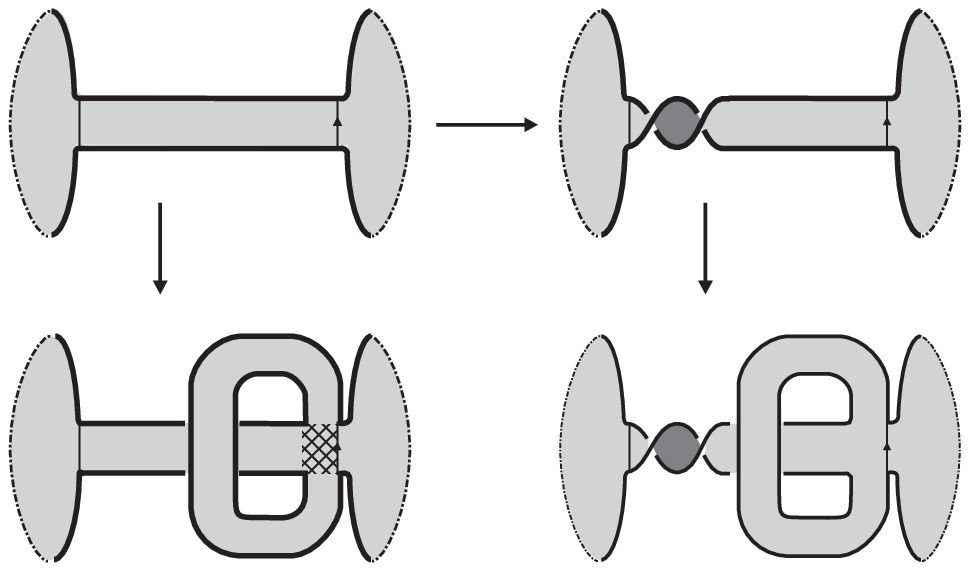}
\put(-4.42,2.03){$\h0v$}\put(-3.7,2.03){$\h1s$}
\put(-3.1,2.03){$\a\quad\h0u$}\put(-0.5,2.03){$\a$}
\put(-2.35,2.2){$T^1_\a$}
\put(-3.7,1.4){$\plumb{\a} A(O,0)$}\put(-1.1,1.4){$\plumb{\a} A(O,0)$}
\put(-2.3,0.5){$\protect\cong$} %
\end{picture}
\caption{$\a$ on $S$ and $T^1_\a(S)$;
$S\plumb{\a} A(O,0)$ and $T^1_\a(S)\plumb{\a} A(O,0)$.
\label{twist across alpha}}
\end{figure}
\end{example}

Of course $S(f_0\connsum f_1,\theta)$ is always a 
connected sum $S(f_0,\theta)\connsum S(f_1,\theta)$.
As the following example shows, for $n>1$, if
$\exp(i\theta)$ and $\Theta$ lie in different components
of $S^1\setminus \crit(f_0\plumb{}f_1)$, then 
$S(f_0\plumb{}f_1,\theta)$ 
need not be a Murasugi sum 
of $S(f_0,\theta)$ and $S(f_1,\theta)$. 

\begin{example}\label{u*u} 
As pictured in Fig.~\ref{uround} (redrawn piecewise-smoothly
in Fig.~\ref{piecewise-smooth U}), and noted in 
Example~\ref{plumbing A(O,0)},
the large Seifert surface of $u$ is $A(O,0)$.
Of course the small Seifert surface of $u$ is two disks.  
By the construction and Prop.~\ref{properties of Murasugi sum of maps},
there is a (moderate, self-indexed) 
$4$--gonal Murasugi sum $u\plumb{}u$ 
with $4$ critical points
whose large Seifert surface $A(O,0)\plumb{} A(O,0)$
is a punctured torus bounded by an unknot; the small 
Seifert surface of $u\plumb{}u$ is a disk, which is 
not a Murasugi sum of two pairs of disks.  
(The situation is pictured 
in Fig.~\ref{piecewise-smooth U*U}.)  
By Example~\ref{plumbing A(O,0)},
there is also a (moderate, self-indexed) 
$4$--gonal Murasugi sum $u\plumb{}o\{2,2\}$ 
with just $2$ critical points (since $o\{2,2\}$ is
a fibration) with the same large and small
Seifert surfaces as $u\plumb{}u$.
\begin{figure}
\centering
\includegraphics[width=4.5in]{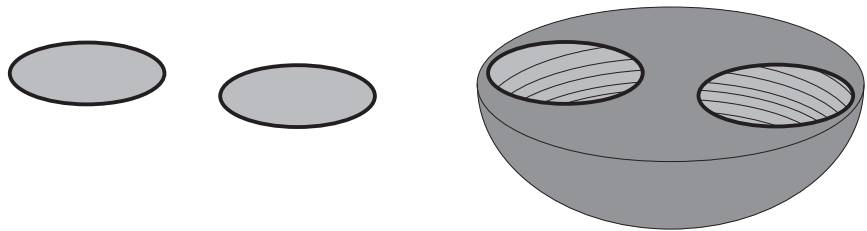}
\caption{Piecewise-smooth Seifert surfaces for $U$ 
(two disks; an annulus)
isotopic to those in Figure \protect\ref{uround}.
\label{piecewise-smooth U}}
\end{figure}
\begin{figure}
\centering
\includegraphics[width=4.5in]{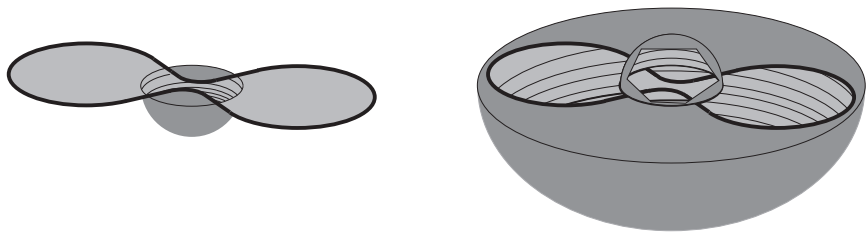}
\caption{Piecewise-smooth Seifert surfaces for 
$U\plumb{}U\protect\cong O$ 
(a disk; a punctured torus, with viewport)
isotopic to a small Seifert surface 
and a large Seifert surface 
of $u\plumb{}u$.
\label{piecewise-smooth U*U}}
\end{figure}
\end{example}

Example~\ref{u*u} also shows that, if no restriction is placed on
the Seifert surfaces along which a Murasugi sum $L_0\plumb{}L_1$ 
is formed, then it can easily happen that 
$\MN(L_0\plumb{}L_1) > \MN(L_0)+\MN(L_1)$.
Much worse is true: according to Hirasawa \cite{Hirasawa:Murasugi},
any knot $K$ bounds a Seifert surface $F$ such that 
$F=F_0\plumb{}F_1$ is a Murasugi sum 
and $K_0\bydefl \Bd F_0$, $K_1\bydefl \Bd F_1$ are unknots; 
since \cite{P-R-W} for every $m$ 
there is a knot $K$ with $\MN(K)>m$, 
$\MN(K_0\plumb{}K_1)-(\MN(K_0)+\MN(K_1))$ can 
be arbitrarily large.
However,
Prop.~\ref{properties of Murasugi sum of maps}
does lead immediately to the following generalization
of the inequality (\ref{subadditivity over connected sum}) 
stated in section~\ref{introduction}.

\begin{cor}\label{subadditivity over Murasugi sum}
If $f_s\from S^3_s\setminus L_s\to S^1$ is a minimal Morse map, 
and $f_0\plumb{(N_0,N_1,\zeta)}f_1$ is 
any Murasugi sum, 
then 
$\MN(L_0\plumb{(N_0,N_1,\zeta)}L_1)\le \MN(L_0)+\MN(L_1)$.
\qed
\end{cor}

\section{Morse maps and free Seifert surfaces}\label{freeness}

For any $X$, let $h_1(X)\bydefl\rank H_1(X;\Z)$.
A Seifert surface $F\subset S^3$ is called 
\emph{free}\label{definition of free}
iff $F$ is connected and $\pi_1(S^3\setminus F)$ 
is a free group; alternatively, $F$ is free iff 
$S^3\setminus \Int \Nb{S^3}{F}$ 
is a handlebody 
$H_g \cong (S^1\times D^2)_1\bdconnsum\dotsm%
{\bdconnsum}(S^1\times D^2)_g$
(necessarily of genus $g=h_1(F)$).
Call the Morse map $f$ \emph{free} if, for every regular value 
$\exp(i\theta)$ of $f$, the Seifert surface $S(f,\theta)$ 
is free.  

\begin{lem}\label{Murasugi sum preserves freeness}
A Murasugi sum $F_0\plumb{} F_1$ of Seifert surfaces $F_0$ 
and $F_1$ is free if and only if $F_0$ and $F_1$ are free.
\end{lem}

\begin{proof} As observed by (for instance) 
Stallings \cite{Stallings}, 
$\pi_1(S^3\setminus F_0\plumb{} F_2)$ is the
free product of $\pi_1(S^3\setminus F_0)$ and   
$\pi_1(S^3\setminus F_1)$.  In particular, 
$\pi_1(S^3\setminus F_0\plumb{} F_1)$ is free
if and only if $\pi_1(S^3\setminus F_0)$ and 
$\pi_1(S^3\setminus F_1)$ are free.
\end{proof} 

It is well known (and obvious) that, if $f$ is a fibration,
then $f$ is free.  
The following proposition is Lemma~4.2 of \cite{P-R-W}.

\begin{prop}
A large Seifert surface of a moderate self-indexed Morse map is free.
\qed
\end{prop}

\begin{cor} 
If $f$ is a moderate self-indexed Morse map
(in particular, if $f$ is a self-indexed minimal Morse map)
then $f$ is free if and only if a small Seifert surface of $f$ is free.
\qed
\end{cor}

Example~\ref{not free}, below, gives a knot $K$ and a 
non-free 
self-indexed minimal Morse map 
$f\from S^3\setminus K\to S^1$. 

Recall that a genus $g$ handlebody $H_g\sub S^3$ 
is \emph{Heegard} if $S^3\setminus\Int H_g$ is 
also a handlebody.  
By a theorem of Waldhausen, any 
two Heegard handlebodies of genus $g$ in $S^3$ are isotopic,
so it is obvious (by considering standard examples)
that, if $H_g\sub S^3$ is a Heegard handlebody 
and $\D\sub\Bd H_g$ is a $2$--disk, then 
$F=\Bd H_g\setminus \Int \D$ is a 
free Seifert surface for the unknot $\Bd\D$.
The converse is also true and easily proved.

\begin{lem}\label{free S for O is completely compressible}
A free Seifert surface $F$ for an unknot 
is of the form $F=\Bd H_g\setminus \Int \D$ 
for some Heegard handlebody $H_g\sub S^3$
and $2$--disk $\D\sub\Bd H_g$. 
\qed
\end{lem}
%


\begin{lem}\label{completely compressible surface is large}
If $S$ is a free Seifert surface for the unknot 
of genus $g$, 
then there is a moderate self-indexed Morse map 
$f\from S^3 \setminus \Bd S\to S^1$ with exactly $h_1(S)=2g$ 
critical points, such that 
$S$ is a large Seifert surface of $f$
and a small Seifert surface of $f$ is a disk.
\end{lem}
\begin{proof} This follows immediately from 
Lemma~\ref{free S for O is completely compressible},
using a suitable modification of the usual dictionary 
between Morse functions (to $\R$) and cobordisms.
Alternatively, and more in the spirit of this paper,
note that (by Prop.~\ref{properties of Murasugi sum of maps}
and Example~\ref{plumbing A(O,0)}) 
a $g$--fold connected sum 
$(u\plumb{}o\{2,2\})\connsum\dotsm\connsum(u\plumb{}o\{2,2\})$
is such a map.
\end{proof}

\begin{prop} If $F$ is a free Seifert surface,
then $F$ is the small Seifert surface of $f$ for 
a free moderate self-indexed Morse
map $f\from S^3 \setminus \Bd F\to S^1$ which has 
exactly $2h_1(F)$ critical points.
\end{prop}

\begin{proof} 
Since $F$ is connected, there exist pairwise disjoint 
proper arcs $\a_s\sub F$, $1\le i\le h_1(F)$, such that 
$(\dotsm((F\cut\a_1)\cut\a_2)\dotsm)\cut  \a_{h_1(F)}$ is a 
$2$-disk. If $T$ denotes the (commutative) composition 
$T^1_{\a_1}\circ T^1_{\a_2}\circ\dotsm\circ T^1_{\a_{h_1(F)}}$,
and $F_1\bydefl T(F)$, then $F_1$ is also free and 
$(\dotsm((F_1\cut\a_1)\cut\a_2)\dotsm)\cut\a_{h_1(F)} %
=(\dotsm((F\cut\a_1)\cut\a_2)\dotsm)\cut  \a_{h_1(F)}$.
The iterated Murasugi sum
\[
F_2\bydefl (\dotsm(F_1\plumb{\a_1}A(O_1,0))\plumb{\a_2}A(O_2,0)\dotsm)%
\plumb{\a_{h_1(F)}}A(O_{h_1(F)},0)
\]
(where each $O_s$ is an unknot) is a Seifert surface
of genus $g(F_2)=2 h_1(F)$, and (as in Example~\ref{plumbing A(O,0)})
$\Bd F_2$ is an unknot. 
By Lemma~\ref{Murasugi sum preserves freeness}, $F_2$ is free.
By Lemma~\ref{free S for O is completely compressible}, 
$F_2=\Bd H_{h_1(F)}\setminus \Int \D$ 
for some Heegard handlebody $H_{h_1(F)}\sub S^3$
and $2$--disk $\D\sub\Bd H_{h_1(F)}$. 
By Prop.~\ref{completely compressible surface is large},
there is a moderate self-indexed Morse map 
$f_2\from S^3 \setminus \Bd F_2\to S^1$ with 
exactly $2g(F_2)=4 h_1(F)$ critical points 
such that $F_2$ is a large Seifert surface of $f$.
But now, 
\[
F_3\bydefl (\dotsm(F_2\plumb{\g(O_1)}A(O'_1,-1))%
\plumb{\g(O_2)}A(O'_2,-1)\dotsm)%
\plumb{\g(O_{h_1(F)})}A(O'_{h_1(F)},-1)
\]
(where each $O'_t$ is an unknot) has boundary isotopic to $\Bd F$.
(The situation in the neighborhood of a single $\a_s$ is
pictured in Fig.~\ref{Hopf plumbing undoes twisting}.)
Since the fibrations $o\{2,2\}$ contribute no new critical
points,
\[
f_3\bydefl (\dotsm(f_2\plumb{}o\{2,2\})%
\plumb{}o\{2,2\}\dotsm)\plumb{}o\{2,2\}
\]
is also moderate and self-indexed with 
exactly $4 h_1(F)$ critical points.  
If the isotopy from $\Bd F_2$ to $\Bd F$ carries
$f_2$ to $f\from S^3 \setminus \Bd F\to S^1$,
then $f$ has the properties required of it.
\begin{figure}
\centering
\includegraphics[width=4.5in]{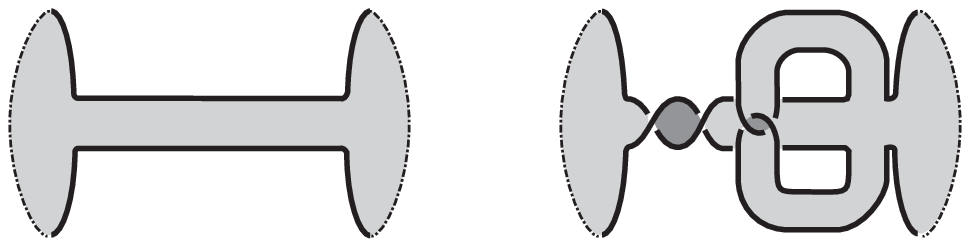}
\caption{$S$ and 
$(T^1_{\a_s}(S)\plumb{{\a_s}}A(O_s,1))\plumb{\g(O_s)}A(O'_s,-1)$
have isotopic boundaries.
\label{Hopf plumbing undoes twisting}}
\end{figure}
\end{proof}

The \emph{free genus} of a knot $K$ is 
$g_f(K)\bydefl \min\{g(F): K=\Bd F$, $F$ is a free Seifert surface$\}$.
More generally, for any link $L$ let 
$h_f(L)\bydefl \min\{h_1(F):L=\Bd F$, $F$ is a free Seifert surface$\}$; 
so, for instance, $h_f(K)=2g_f(K)$ if $L=K$ is a knot, 
and $h_f(L)=1$ iff $L=\Bd A(O,n)$ for some $n\in\Z$, $n\ne 0$.

\begin{cor}\label{MN bounded by four times the free genus} 
For any link $L$, $\MN(L)\le 2h_f(L)$.  In particular, 
for any knot $K$, $\MN(K)\le 4g_f(K)$.
\qed 
\end{cor}

Let $g(K)$, as usual, denote the genus of a knot $K$.
Of course $g_f(K)\ge g(K)$ for every knot $K$; 
there are knots $K$ with $g(K)=1$ and $g_f(K)$ arbitrarily 
large \cite{Moriah:free,Livingston:free}.
For many knots $K$ (satisfying a suitable condition on
the Alexander polynomial $\Delta_K(t)$), it follows 
from \cite{P-R-W} that $\MN(K)\ge 2g(K)$.  I know no 
example of a knot $K$ for which it can be shown that 
$\MN(K)>2g(K)$ (the techniques of \cite{P-R-W}, as they 
stand, cannot be used to do this).  Some examples follow 
in which $\MN(K)=2g(K)$.

\begin{lem}\label{MN(A(O,n))}
$\MN(\Bd A(O,n))=2$ for all $n\ne\pm1$.
\end{lem}
\begin{proof} This generalization of Prop.~\ref{MN(U)=2}
is implicit in Example~\ref{plumbing A(O,0)}.
\end{proof}

\begin{example}\label{MN of twist knots}
A twist knot is, by definition, the boundary of a 
Seifert surface of genus $1$ of the form 
$A(O,n)\plumb{}A(O,\mp1)$.  
By Cor.~\ref{subadditivity over Murasugi sum},
Lem.~\ref{MN(A(O,n))}, and the fiberedness of 
$A(O,\mp1)$, it follows that, if $K$ is a twist knot, 
then $\MN(K)=2$ unless $K$ is fibered,
and in every case $K$ has a free minimal Morse function.
\end{example}

\begin{example}\label{not free}
Let $S$ be a surface of genus $1$ with connected boundary.
If $\a\sub S$ is a non-separating (i.e., boundary-incompressible)
proper arc, then $S\cut\a$ is an annulus; conversely,
if $c\sub S$ is a non-separating simple closed curve,
then there is a non-separating proper arc $\a\sub S$
such that $S\cut\a=\Nb{S}{c}$.
By judicious choices of $\a$ on the (genus $1$) 
fiber surface $F$ of one of the non-trivial fibered 
twist knots (that is, either a positive or negative 
trefoil knot $O\{2,\pm3\}=\Bd (A(O,\mp1)\plumb{} A(O,\mp1)$
or the figure-8 knot $\Bd (A(O,1)\plumb{} A(O,-1)$),
one may find annuli $A(K,n)\sub F$ for knots $K$
of infinitely many isotopy types.  
Any such annulus $A(K,n)$ is the small
surface, and $(T^{\pm1}_\a(F)\plumb{\a}A(O,0))\plumb{\g(O)}A(O',\mp1)$
is the large surface, of a moderate Morse map 
$f\from S^3\setminus\Bd A(K,n)\to S^1$ with 
exactly $2$ critical points.  Such a Morse map is 
minimal and not free.
By plumbing on another copy of $A(O,\mp1)$
we obtain infinitely many doubled knots
$D(K,n,{\mp1})\bydefl \Bd(A(K,n)\plumb{}A(O,\mp1))$
with $\MN(D(K,n,{\mp1}))=2$, each of which 
has a non-free minimal Morse map.
\begin{figure}
\centering
\includegraphics[width=3.5in]{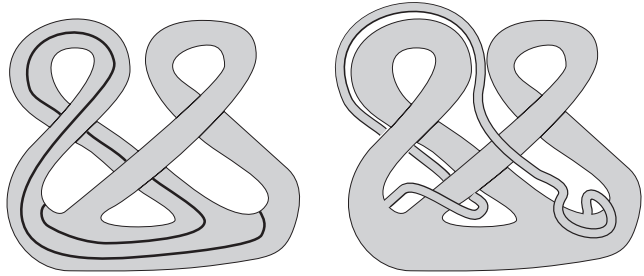}
\caption{A proper arc $\a$ on $S(o\{2,3\},\theta_0)$;
a free surface $S(o\{2,3\},\theta_0)\plumb{\a}A(O,0)$
isotopic to a large Seifert surface 
$S(o\{2,3\}\plumb{}u,\theta_0)$.\label{unfree}}
\end{figure}
\begin{figure}
\centering
\includegraphics[width=3.5in]{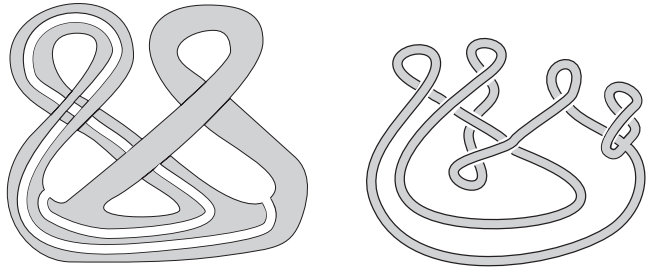}
\caption{A non-free knotted annulus
$A(O\{2,3\},-1)$ isotopic to a small Seifert 
surface $S(o\{2,3\}\plumb{}u,\theta_1)$; 
a non-free punctured torus 
$A(O\{2,3\},-1)\plumb{}A(O,-1)$
which is a large Seifert surface of 
a minimal Morse function 
$(o\{2,3\}\plumb{}u)\plumb{}o\{2,2\}$
for the doubled 
knot $D(O\{2,3\},-1,{+})$.
\label{unfree2}}
\end{figure}
\end{example}

\begin{hist}
Neuwirth \cite{Neuwirth:algebraic} was 
apparently the first (1960) to consider, in effect, the
notion of ``free Seifert surface'': in a footnote, he 
called a Seifert surface $F$ ``algebraically knotted'' 
if $\pi_1( S^3\setminus F)$ is not free.  
Neuwirth's language was adapted by 
Murasugi \cite{Murasugi:commutator}, 
Lyon \cite{Lyon:1972}, and others.
By 1976, Problem~1.20 (attributed to Giffen and Siebenmann)
in Kirby's problem list \cite{Kirby:problems}
uses the phrase ``free Seifert surface'', 
points out that $F$ is free iff $S^3\setminus F$ is an open handlebody,
and apparently introduces the terminology ``free genus''.  
A number of later authors
(e.g., Moriah \cite{Moriah:free},
Livingston \cite{Livingston:free},
M.~Kobayashi and T.~Kobayashi \cite{Kobayashi-Kobayashi})
have studied free Seifert surfaces.
\end{hist}

\section{The inhomogeneity and Morse-Novikov number of a closed braid}
\label{inhomogeneity}

Let $\s_1,\dots,\s_{n-1}$ be the standard generators of the
$n$--string braid group $B_n$.  A \emph{braidword} of length
$k$ in $B_n$ is a $k$--tuple $\bword=(b(1),\dots,b(k))$
such that $b(s)=\s_{i_\bword(s)}^{\e_\bword(s)}$ for some 
$(i_\bword,\e_\bword)\from \{1,\dots,k\}\to\{1,\dots,n-1\}\times\{-1,1\}$.  
The \emph{braid} of $\bword$ is $\b(\bword)\bydefl b(1)\dotsm\b(k)$.
The \emph{closed braid} of $\bword$ is 
$\widehat\b(\bword)\bydefl \widehat{\b(\bword)}$, where
the closure of a braid $B_n$ is (as usual) a certain
link in $S^3$.
A braidword $\bword$ for which $i_\bword$ is surjective
is \emph{strict}.  If $\bword$ is not strict, then 
$\widehat\b(\bword)$ is a split link
(i.e., $\pi_2(S^3\setminus \widehat\b(\bword)\ne\{0\}$).

\begin{definition*}\label{definition of inhomogeneity}
Let $\bword$ be a strict braidword in $B_n$.  For 
$i\in \{1,\dots,n-1\}$ and $\e\in\{-1,1\}$, let
$v(i,\e,\bword)\bydefl \card{(i_\bword,\e_\bword)^{-1}(i,\e)}$.
The \emph{inhomogeneity of $\bword$}
is 
\[
I(\bword)\bydefl \sum_{i=1}^{n-1} \min\{v(i,-1,\bword),v(i,1,\bword)\}
\]
and the \emph{inhomogeneity of $\b\in B_n$}
is $I(\b)\bydefl \min\{I(\bword):\b=\b(\bword), \bword$ is strict$\}$.
A braid with inhomogeneity $0$ is \emph{homogeneous}.
\end{definition*}

The following result generalizes 
Theorem~2 of \cite{Stallings}, which states that the 
closure of a homogeneous braid is a fibered link.

\begin{prop}\label{MN bounded by inhomogeneity}
For all $n$, for all $\b\in B_n$, $\MN(\widehat\b)\le 2I(\b)$.
\end{prop}
\begin{proof} There is nothing to prove if $n=1$.
Let $n=2$.  Given a strict braidword 
$\bword=(\s_1^{\e_\bword(1)},\dots,\s_1^{\e_\bword(k)})$ 
in $B_2$, application of Seifert's algorithm to the
corresponding closed braid diagram 
for $\widehat\beta(\bword)$ 
produces a Seifert surface $S(\bword)$ for $\widehat\beta(\bword)$ 
naturally equipped with a $(0,1)$--handle decomposition
$S=(\h01\cup\h02)\cup \bigcup_{t=1}^k \h1t$
such that $\h01$ and $\h02$ are $2$--disks embedded 
in parallel planes in $\R^3\sub S^3$ and $\h1t$ joins
$\h01$ to $\h02$ with a single half-twist of sign $\e_\bword(t)$.
In case $\bword$ is homogeneous, $S(\bword)$ is (as is well
known) a fiber surface.  If $\bword$ is not homogeneous, then
(up to a cyclic permutation of $\bword$, which changes neither
$\b(\bword)$ nor the isotopy type of $S(\bword)$) there exists
$t<k$ such that $\e_\bword(t)=1=-\e_\bword(t+1)$, and 
$S(\bword)$ is isotopic to an annulus plumbing 
$S(\bword')\plumb{}A(O,0)$, where 
$\bword'\bydefl(b(1),\dots,b(t-1),b(t+2),\dots,b(k))$.
(The situation is pictured in 
Fig.~\ref{unplumbing an inhomogeneous braided surface of degree 2}.)  
By Prop.~\ref{MN(U)=2}, 
Cor.~\ref{subadditivity over Murasugi sum}, 
and induction on $I(\b)$, 
the proposition is true for $n=2$.
Finally, as observed by Stallings \cite{Stallings}, 
for all $n\ge 2$, application of Seifert's algorithm 
to the braidword diagram of the closed braid of 
a strict braidword in $B_n$ produces a Seifert surface 
which is an iterated Murasugi sum of $n-1$ Seifert surfaces 
of strict braidwords in $B_2$.
By Cor.~\ref{subadditivity over Murasugi sum}
and induction, the proposition is true for all $n$.
\begin{figure}
\centering
\includegraphics[width=3.5in]{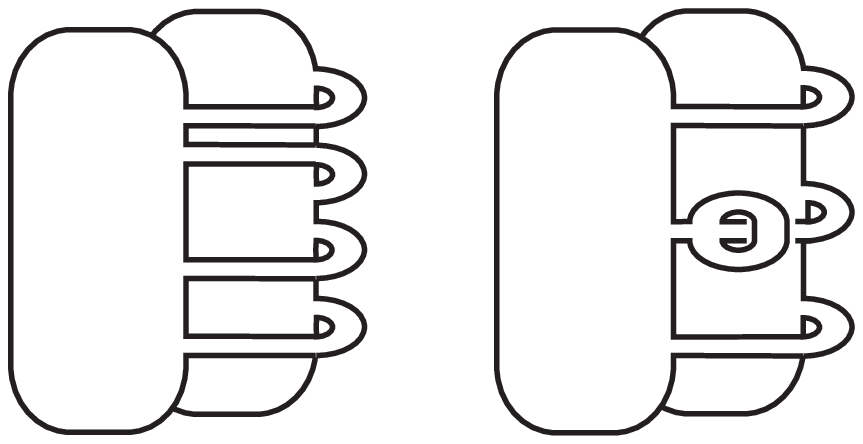}
\caption{$S(\s_1,\s_1,\s_1^{-1},\s_1)$ is isotopic 
by handle-sliding to
$S(\s_1,\s_1,\s_1)\plumb{}A(O,0)$.
\label{unplumbing an inhomogeneous braided surface of degree 2}}
\end{figure}
\end{proof}

Just like Theorem~2 of \cite{Stallings}, which 
can be generalized to closed $\Ts$-homogeneous braids
for any espalier $\Ts$ on vertices $\{1,\dots,n\}$
(see \cite{espaliers}), Prop.~\ref{MN bounded by inhomogeneity}
can be generalized in terms of a suitably defined 
notion of the ``$\Ts$-inhomogeneity'' $I_\Ts(\b)$
of a braid $\b\in B_n$.  Similarly, the explicit construction
(by lifting the fibration $o$ through a suitable
branched cover $S^3\to S^3$) of a fibration 
$f\from S^3\setminus\widehat\b\to S^1$, 
given in \cite{someknot} for a closed homogeneous braid 
$\widehat\b$ and easily generalized (as mentioned in
\cite{specpos}) to closed $\Ts$-homogeneous braids,
actually gives for any $\Ts$-bandword $\bword$ (in 
particular for an ordinary bandword) an explicit
self-indexed moderate Morse map 
$f\from S^3\setminus\widehat\b(\bword)\to S^1$
with $2I_\Ts(\bword)$ critical points.

\providecommand{\bysame}{\leavevmode\hbox to3em{\hrulefill}\thinspace}

\end{document}